\pdfoutput=1
\documentclass[11pt,a4paper]{article}

\usepackage[left=30mm,right=29mm,top=29mm,bottom=31mm]{geometry}
\usepackage[T1]{fontenc}
\usepackage{lmodern}
\usepackage{amsmath,amssymb,amsthm}
\usepackage{microtype}
\usepackage[hyperfootnotes=false]{hyperref}
\usepackage{enumitem}
\usepackage{needspace}
\allowdisplaybreaks[2]

\title{\Large On a question of A.N. Skiba concerning the modularity\\ of the lattice of $\sigma$-permutable subgroups}
\author{\normalsize Chenxiao Lin\\
\small School of Mathematics and Statistics, Jiangsu Normal University,\\
\small Xuzhou 221116, Jiangsu, China\\
\small E-mail: 3020243174@jsnu.edu.cn}
\date{}
\hypersetup{pdfborder={0 0 1},linkbordercolor={0 1 0},citebordercolor={0 1 0},urlbordercolor={0 1 1},pdftitle={On a question of A.N. Skiba concerning the modularity of the lattice of sigma-permutable subgroups},pdfauthor={Chenxiao Lin},pdfsubject={Finite groups; sigma-permutable subgroups; subgroup lattices},pdfkeywords={finite group; Hall subgroup; sigma-permutable subgroup; subgroup lattice; modular lattice}}

\theoremstyle{plain}
\newtheorem{theorem}{Theorem}

\newtheorem{lemma}{Lemma}[section]
\newtheorem{proposition}[lemma]{Proposition}
\newtheorem{corollary}[lemma]{Corollary}
\newtheorem{question}[lemma]{Question}
\theoremstyle{definition}
\newtheorem{example}[lemma]{Example}

\setlist[enumerate]{label=\textup{(\roman*)},leftmargin=2.2em,itemsep=2pt,topsep=4pt}

\newcommand{\Lall}{\mathcal L}
\newcommand{\Lper}{\mathcal L_{\sigma\mathrm{per}}}
\newcommand{\Lsi}{\mathcal L_{\sigma_i}}
\newcommand{\gen}[1]{\langle #1\rangle}

\begin{document}
\maketitle
\begin{abstract}
Let $\sigma=\{\sigma_i\mid i\in I\}$ be some partition of the set of all primes and let $G$ be a finite group. A subgroup $A$ of $G$ is said to be $\sigma$-permutable in $G$ if $G$ has a Hall $\sigma_i$-subgroup for every $i\in I$ and $A$ permutes with every Hall $\sigma_i$-subgroup of $G$ for all $i\in I$. We describe the finite groups $G$ for which the lattice $\Lper(G)$ of all $\sigma$-permutable subgroups of $G$ is modular. This answers the first part of a question of A.N. Skiba. We also prove that if every Hall $\sigma_i$-subgroup of $G$ is nilpotent, then $\Lper(G)$ is modular if and only if every two members of $\Lper(G)$ are permutable.
\end{abstract}

{\let\thefootnote\relax
}

\section{Introduction}

Throughout this paper, all groups are finite and $G$ always denotes a finite group. Moreover, $\mathbb P$ is the set of all primes, $\pi(G)$ is the set of all primes dividing $|G|$, and $\Lall(G)$ is the lattice of all subgroups of $G$. The subgroups $A$ and $B$ of $G$ are said to be \emph{permutable} if $AB=BA$; in this case we also say that $A$ permutes with $B$.

In what follows, $\sigma=\{\sigma_i\mid i\in I\}$ is some partition of $\mathbb P$, that is, $\mathbb P=\bigcup_{i\in I}\sigma_i$ and $\sigma_i\cap\sigma_j=\varnothing$ for all $i\ne j$, and $\sigma(G)=\{\sigma_i\mid\sigma_i\cap\pi(G)\ne\varnothing\}$. The group $G$ is said to be \emph{$\sigma$-full} if $G$ has a Hall $\sigma_i$-subgroup for every $i\in I$. A subgroup $A$ of $G$ is said to be \emph{$\sigma$-permutable} in $G$ \cite{Skiba2015} if $G$ is $\sigma$-full and $A$ permutes with every Hall $\sigma_i$-subgroup of $G$ for all $i\in I$. We write $\Lper(G)$ to denote the set of all $\sigma$-permutable subgroups of $G$. If $G$ is $\sigma$-full, then $\Lper(G)$ is a sublattice of $\Lall(G)$ \cite[Theorem~A]{HuHuangSkiba} (see also \cite[Theorem~B]{SkibaArxiv}), so that the meet and the join in $\Lper(G)$ are the intersection and the subgroup generated, respectively.

\begin{question}[{See Question~1.4 in \cite{HuHuangSkiba}; see also \cite[Question~1.3]{SkibaArxiv} and \cite[Question~6]{SkibaNote}}]\label{q:skiba}
Under what conditions the lattice $\Lper(G)$ is modular or distributive?
\end{question}

The second part of Question~\ref{q:skiba} was answered in \cite[Theorem~B]{HuHuangSkiba} (see also \cite[Theorem~A]{SkibaArxiv}), where the $\sigma$-full groups $G$ with distributive lattice $\Lper(G)$ were described. On the other hand, it was noted in \cite{HuHuangSkiba} that the conditions under which the lattice $\Lper(G)$ is modular were not known. In this paper we give an answer to the first part of Question~\ref{q:skiba}.

In the classical case when $\sigma=\{\{2\},\{3\},\{5\},\ldots\}$, a subgroup of $G$ is $\sigma$-permutable if and only if it is $S$-permutable, that is, it permutes with all Sylow subgroups of $G$. In this case, Skiba \cite{SkibaArch} proved that the lattice $\Lall_S(G)$ of all $S$-permutable subgroups of $G$ is modular if and only if every two members of $\Lall_S(G)$ are permutable (see also \cite[Theorem~2.5]{SelkinSkiba}). Earlier, Kimber \cite{Kimber} proved a similar result for the lattice of all $\Sigma$-permutable subgroups of a soluble group with a Hall system $\Sigma$ (see also \cite[Theorem~27]{SkibaNote}). In general, however, the modularity of $\Lper(G)$ does not imply that its members are permutable. For example, if $\pi(G)\subseteq\sigma_i$ for some $i$, then $\Lper(G)=\Lall(G)$; in particular, the lattice $\Lall(S_3)$ is modular, but no two distinct subgroups of order $2$ of $S_3$ are permutable.

Let $Q$ be a $\sigma$-full group and $i\in I$. We write $O_{\sigma_i}(Q)$ to denote the largest normal $\sigma_i$-subgroup of $Q$, and $O^{\sigma_i}(Q)$ to denote the intersection of all normal subgroups $K$ of $Q$ such that $Q/K$ is a $\sigma_i$-group. Let
\begin{equation}\label{eq:operators}
P=O_{\sigma_i}(Q)\qquad\text{and}\qquad D=O^{\sigma_i}(Q).
\end{equation}
For every normal subgroup $R$ of $D$ with $R\le P\cap D$, we put
\begin{equation}\label{eq:section}
P_R=\{x\in P\mid [x,D]\le R\},
\end{equation}
where $[x,D]$ is the subgroup generated by all commutators $[x,d]=x^{-1}d^{-1}xd$, $d\in D$. Then $P_R$ is a subgroup of $P$ and $R$ is normal in $P_R$ (see Proposition~\ref{prop:centralizer} below). In fact, $P_R/R$ is the set of all elements of $N_P(R)/R$ which centralize $D/R$ in $N_Q(R)/R$.

\begin{theorem}\label{thm:main}
Let $G$ be a $\sigma$-full group. Then the lattice $\Lper(G)$ is modular if and only if for every normal subgroup $N$ of $G$, every $i\in I$ and every normal subgroup $R$ of $D$ with $R\le P\cap D$, where $P$ and $D$ are defined by~\eqref{eq:operators} for $Q=G/N$, the group $P_R/R$ is a direct product of groups of pairwise coprime orders each of which is of one of the following types:
\begin{enumerate}
\item $Q_8\times E$, where $E$ is an elementary abelian $2$-group;
\item a $p$-group $X=U\gen{t}$, where $U$ is an abelian normal subgroup of $X$ and $t^{-1}ut=u^{1+p^s}$ for all $u\in U$ and some integer $s\ge1$, with $s\ge2$ if $p=2$;
\item a semidirect product $(C_p)^n\rtimes C_{q^k}$, where $n,k\ge1$, $p$ and $q$ are distinct primes, and a generator of $C_{q^k}$ acts on $(C_p)^n$ as the multiplication by an element of order $q$ of $\mathbb F_p^\times$.
\end{enumerate}
\end{theorem}

Here the identity group is regarded as the empty direct product. In~(ii), the subgroups $U$ and $\gen{t}$ may have nontrivial intersection, and in~(iii) we have $q\mid p-1$. The groups of types~(i)--(iii) are exactly the factors in Iwasawa's description of the finite groups with modular subgroup lattice \cite[Theorems~2.3.1 and~2.4.4]{Schmidt}. Thus Theorem~\ref{thm:main} states that $\Lper(G)$ is modular if and only if the subgroup lattice of every group $P_R/R$ in Theorem~\ref{thm:main} is modular.

\begin{theorem}\label{thm:nilpotent}
Let $G$ be a $\sigma$-full group. Suppose that $O_{\sigma_i}(G/N)$ is nilpotent for every normal subgroup $N$ of $G$ and every $i\in I$. Then the lattice $\Lper(G)$ is modular if and only if every two members of $\Lper(G)$ are permutable. In particular, this is true if every Hall $\sigma_i$-subgroup of $G$ is nilpotent for all $i\in I$.
\end{theorem}

In the classical case when $\sigma=\{\{2\},\{3\},\{5\},\ldots\}$, every Hall $\sigma_i$-subgroup of $G$ is a Sylow subgroup, so Theorem~\ref{thm:nilpotent} contains the above result of Skiba.

\begin{corollary}[{Skiba \cite{SkibaArch}}]\label{cor:sylow}
The lattice $\Lall_S(G)$ is modular if and only if every two members of $\Lall_S(G)$ are permutable.
\end{corollary}

Example~\ref{ex:dihedral} shows that the nilpotency hypothesis in Theorem~\ref{thm:nilpotent} cannot be omitted.

In Section~2 we reduce the modularity of $\Lper(G)$ to the modularity of the lattices of all $\sigma$-permutable $\sigma_i$-subgroups of the quotients of $G$ (Proposition~\ref{prop:quotient}). In Section~3 we prove Theorem~\ref{thm:main}, and in Section~4 we prove Theorem~\ref{thm:nilpotent} and give a sufficient condition for the modularity of $\Lper(G)$ in terms of a normal Hall subgroup of $G$. In Section~5 we show that the conditions in Proposition~\ref{prop:quotient} and Theorem~\ref{thm:main} cannot be weakened.

All unexplained notation and terminology are standard. The reader is referred to \cite{Schmidt,DoerkHawkes} if necessary.

\section{Preliminaries}

In what follows, $A_G$ denotes the core of a subgroup $A$ in $G$. Recall that $G$ is said to be \emph{$\sigma$-nilpotent} if $G$ is the direct product of its Hall $\sigma_i$-subgroups, $i\in I$. If $G$ is $\sigma$-nilpotent and $H_i$ is the Hall $\sigma_i$-subgroup of $G$, then every subgroup $A$ of $G$ is $\sigma$-nilpotent and $A=\prod_{i\in I}(A\cap H_i)$, since every element of $G$ is the product of its $\sigma_i$-components, which are powers of this element. Note also that if $G$ is $\sigma$-full and $N$ is a normal subgroup of $G$, then $G/N$ is $\sigma$-full, since the image of every Hall $\sigma_i$-subgroup of $G$ in $G/N$ is a Hall $\sigma_i$-subgroup of $G/N$.

Let $P$ and $D$ be normal subgroups of a group $Q$. We write $\Lall_D(P)$ to denote the set of all $D$-invariant subgroups of $P$; clearly, $\Lall_D(P)$ is a sublattice of $\Lall(P)$. If $G$ is $\sigma$-full and $i\in I$, then $\Lsi(G)$ denotes the set of all $\sigma$-permutable $\sigma_i$-subgroups of $G$ (see \cite{SkibaArxiv}).

\begin{lemma}\label{lem:skiba}
Let $G$ be a $\sigma$-full group and $N$ a normal subgroup of $G$.

\textup{(1)} \textup{(See \cite[Theorem~A]{HuHuangSkiba} or \cite[Theorem~B]{SkibaArxiv}.)} $\Lper(G)$ is a sublattice of $\Lall(G)$.

\textup{(2)} \textup{(See \cite[Theorem~B]{Skiba2015} or \cite[Lemma~2.2]{HuHuangSkiba}.)} If $A\in\Lper(G)$, then $A/A_G$ is $\sigma$-nilpotent.

\textup{(3)} \textup{(See \cite[Corollary~2.4]{HuHuangSkiba} or \cite[Corollary~2.6]{SkibaArxiv}.)} A $\sigma$-nilpotent subgroup $A$ of $G$ belongs to $\Lper(G)$ if and only if every Hall $\sigma_i$-subgroup of $A$ belongs to $\Lper(G)$ for all $i\in I$.

\textup{(4)} \textup{(See \cite[Lemma~3.2]{SkibaArxiv} and \cite[Lemma~3.1]{HuHuangSkiba}.)} If $A\in\Lper(G)$, then $AN/N\in\Lper(G/N)$. If $N\le V\le G$ and $V/N\in\Lper(G/N)$, then $V\in\Lper(G)$. Hence the lattice $\Lper(G/N)$ is isomorphic to the interval of $\Lper(G)$ consisting of all members containing $N$.
\end{lemma}

\begin{lemma}\label{lem:block}
Let $G$ be a $\sigma$-full group, $i\in I$, $P=O_{\sigma_i}(G)$ and $D=O^{\sigma_i}(G)$. Then $\Lsi(G)=\Lall_D(P)$. In particular, $\Lsi(G)$ is a sublattice of $\Lper(G)$ \textup{(}see also \cite[Proposition~3.3]{HuHuangSkiba}\textup{)}.
\end{lemma}

\begin{proof}
First we show that $P$ is the intersection of all Hall $\sigma_i$-subgroups of $G$. Indeed, this intersection is a normal $\sigma_i$-subgroup of $G$, so it is contained in $P$. On the other hand, if $H$ is a Hall $\sigma_i$-subgroup of $G$, then $HP$ is a $\sigma_i$-subgroup of $G$ containing $H$, so $HP=H$ and $P\le H$.

Next we show that $D$ is generated by all Hall $\sigma_j$-subgroups of $G$ with $j\ne i$. Let $X$ be the subgroup generated by these subgroups. Then $X$ is normal in $G$ and $|G:X|$ is a $\sigma_i$-number, so $D\le X$. On the other hand, if $K$ is a normal subgroup of $G$ such that $G/K$ is a $\sigma_i$-group, then every Hall $\sigma_j$-subgroup of $G$ with $j\ne i$ is contained in $K$, so $X\le K$. Hence $X\le D$, and so $X=D$.

Now let $A\in\Lsi(G)$. If $H$ is a Hall $\sigma_i$-subgroup of $G$, then $AH$ is a $\sigma_i$-subgroup of $G$ containing $H$, so $A\le AH=H$. Hence $A\le P$. Let $K$ be a Hall $\sigma_j$-subgroup of $G$, where $j\ne i$. Then $AK$ is a subgroup of $G$ and $P\cap AK=A(P\cap K)=A$, since $A\le P$ and $P\cap K=1$. Since $P$ is normal in $G$, it follows that $A$ is normal in $AK$, that is, $K$ normalizes $A$. Therefore $D$ normalizes $A$, and so $A\in\Lall_D(P)$.

Conversely, let $A\in\Lall_D(P)$. Then every Hall $\sigma_i$-subgroup $H$ of $G$ contains $P$, so $AH=H=HA$. Moreover, every Hall $\sigma_j$-subgroup $K$ of $G$ with $j\ne i$ is contained in $D$ and so normalizes $A$, which implies that $AK=KA$. Hence $A\in\Lsi(G)$.

Finally, $\Lall_D(P)$ is closed under taking intersections and joins, so $\Lsi(G)$ is a sublattice of $\Lper(G)$. The lemma is proved.
\end{proof}

\begin{lemma}\label{lem:lift}
Let $G$ be a $\sigma$-full group and $N$ a normal subgroup of $G$ such that the lattice $\Lper(G/N)$ is modular. Let $A,B,C\in\Lper(G)$ with $A\le C$. If either $N\le A$ or $N\le B$, then
\begin{equation}\label{eq:modularity}
\gen{A,B}\cap C=\gen{A,B\cap C}.
\end{equation}
\end{lemma}

\begin{proof}
Let $T=\gen{A,B}$, $E=T\cap C$ and $F=\gen{A,B\cap C}$. Then $F\le E$, and the images of $A$, $B$ and $C$ in $G/N$ belong to $\Lper(G/N)$ by Lemma~\ref{lem:skiba}(4).

First assume that $N\le A$. Then $N\le F\le E$ and $BN\cap C=(B\cap C)N$ by the Dedekind identity, since $N\le C$. Hence the modular law in $\Lper(G/N)$, applied to $A/N\le C/N$ and $BN/N$, gives $E/N=F/N$, and so $E=F$.

Now assume that $N\le B$. Then $N\le T$, and so $T\cap CN=(T\cap C)N=EN$ and $B\cap CN=(B\cap C)N$ by the Dedekind identity. Hence the modular law in $\Lper(G/N)$, applied to $AN/N\le CN/N$ and $B/N$, gives $EN=FN$. Since $C\cap N\le B\cap C\le F\le E$ and $E\le C$, the Dedekind identity implies that
\[
E=E(C\cap N)=C\cap EN=C\cap FN=F(C\cap N)=F.
\]
The lemma is proved.
\end{proof}

\begin{proposition}\label{prop:quotient}
Let $G$ be a $\sigma$-full group. Then the lattice $\Lper(G)$ is modular if and only if the lattice $\Lsi(G/N)$ is modular for every normal subgroup $N$ of $G$ and every $i\in I$.
\end{proposition}

\begin{proof}
If $\Lper(G)$ is modular, then $\Lper(G/N)$ is modular for every normal subgroup $N$ of $G$ by Lemma~\ref{lem:skiba}(4), and so its sublattice $\Lsi(G/N)$ is modular by Lemma~\ref{lem:block}.

Now assume that $\Lsi(G/N)$ is modular for all $N$ and $i$, but $\Lper(G)$ is not modular, and let $G$ be a counterexample of minimal order. Since every quotient of $G/N$ is a quotient of $G$, the choice of $G$ implies that $\Lper(G/N)$ is modular for every normal subgroup $N\ne1$ of $G$. Let $A,B,C\in\Lper(G)$ be such that $A\le C$ and~\eqref{eq:modularity} does not hold. If $A_G\ne1$, then Lemma~\ref{lem:lift} with $N=A_G$ implies that~\eqref{eq:modularity} holds, a contradiction. Hence $A_G=1$, and similarly $B_G=1$. Therefore $A$ and $B$ are $\sigma$-nilpotent by Lemma~\ref{lem:skiba}(2).

Let $A_i$ and $B_i$ be the Hall $\sigma_i$-subgroups of $A$ and $B$, respectively. Then $A_i,B_i\in\Lsi(G)$ by Lemma~\ref{lem:skiba}(3), and so $A_i,B_i\le O_{\sigma_i}(G)$ by Lemma~\ref{lem:block}. Since $O_{\sigma_i}(G)\cap O_{\sigma_j}(G)=1$ for all $i\ne j$, these normal subgroups centralize each other. Hence
\[
T=\gen{A,B}=\prod_{i\in I}T_i,\qquad\text{where }T_i=\gen{A_i,B_i}\le O_{\sigma_i}(G),
\]
and only finitely many of the subgroups $T_i$ are nontrivial. Let $C_0=C\cap T$. Then $C_0\in\Lper(G)$ by Lemma~\ref{lem:skiba}(1), and $C_0=\prod_{i\in I}C_i$, where $C_i=C_0\cap T_i$ is the Hall $\sigma_i$-subgroup of $C_0$. Moreover, $C_i\in\Lsi(G)$ by Lemma~\ref{lem:skiba}(3). Since $A_i\le C_i$ and the lattice $\Lsi(G)$ is modular, we have
\[
C_i=T_i\cap C_i=\gen{A_i,B_i\cap C_i}.
\]
Therefore $C\cap T=C_0=\gen{A,B\cap C_0}=\gen{A,B\cap C}$, a contradiction. The proposition is proved.
\end{proof}

\section{\texorpdfstring{Proof of Theorem~\ref{thm:main}}{Proof of Theorem A}}

\begin{proposition}\label{prop:centralizer}
Let $P$ and $D$ be normal subgroups of a group $Q$, and let $R$ be a normal subgroup of $D$ with $R\le P\cap D$. Define $P_R$ by~\eqref{eq:section}. Then $P_R$ is a $D$-invariant subgroup of $P$ and $R$ is normal in $P_R$. Moreover, the lattice $\Lall_D(P)$ is modular if and only if the lattice $\Lall(P_R/R)$ is modular for every such $R$.
\end{proposition}

\begin{proof}
Let $x\in P_R$ and $d\in D$. Then $x^{-1}dx=d[x,d]^{-1}\in dR$. In particular, if $d\in R$, then $x^{-1}dx\in R$, so $x$ normalizes $R$. Hence $P_R$ is the kernel of the action of $N_P(R)$ on $D/R$ by conjugation, and so $P_R$ is a subgroup of $P$ containing $R$ and normalizing $R$. Moreover, $P_R$ is $D$-invariant, since $D$ normalizes both $P$ and $R$. If $R\le U\le P_R$, then $U$ is $D$-invariant, since $u^d=u[u,d]\in U$ for all $u\in U$ and $d\in D$. Therefore the lattice $\Lall(P_R/R)$ is isomorphic to the interval $[R,P_R]$ of $\Lall_D(P)$. Hence if $\Lall_D(P)$ is modular, then so is $\Lall(P_R/R)$.

Conversely, assume that the lattices $\Lall(P_R/R)$ are modular for all $R$. Let $A,B,C\in\Lall_D(P)$ with $A\le C$, and let
\[
T=\gen{A,B},\quad E=T\cap C,\quad F=\gen{A,B\cap C},\quad a=A\cap D,\quad b=B\cap D\quad\text{and}\quad R=ab.
\]
Then $a$ and $b$ are normal in $D$, so $R$ is a normal subgroup of $D$ with $R\le P\cap D$. Since $D$ is normal in $Q$ and normalizes $A$ and $B$, we have $[A,D]\le a$ and $[B,D]\le b$. Hence $A,B\le P_R$, so $R\le T\le P_R$ and $R$ is normal in $T$. In particular, the lattice $\Lall(T/R)$ is modular.

Since $a\le A\le E$ and $b\le B$, the Dedekind identity gives
\begin{equation}\label{eq:section-meet}
BR\cap ER=Ba\cap Eb=(B\cap Eb)a=(B\cap E)ba=(B\cap E)R.
\end{equation}
Now $B\cap E=B\cap C$ and $A\le E$. Hence the modular law in $\Lall(T/R)$, applied to $AR/R\le ER/R$ and $BR/R$, together with~\eqref{eq:section-meet} gives $ER=FR$. Moreover, $E\cap R=a(E\cap b)\le F$, since $a\le A$ and $E\cap b\le B\cap C$. Therefore
\[
E=E\cap FR=F(E\cap R)=F
\]
by the Dedekind identity. Thus the modular law holds in $\Lall_D(P)$. The proposition is proved.
\end{proof}

\begin{proof}[Proof of Theorem~\ref{thm:main}]
By Proposition~\ref{prop:quotient}, the lattice $\Lper(G)$ is modular if and only if the lattice $\Lsi(Q)$ is modular for every quotient $Q=G/N$ of $G$ and every $i\in I$. By Lemma~\ref{lem:block}, we have $\Lsi(Q)=\Lall_D(P)$, where $P$ and $D$ are defined by~\eqref{eq:operators}. Hence, by Proposition~\ref{prop:centralizer}, $\Lper(G)$ is modular if and only if the lattice $\Lall(P_R/R)$ is modular for all $N$, $i$ and $R$ as in the theorem.

Recall that, by Iwasawa's theorem (see \cite[Theorems~2.3.1 and~2.4.4]{Schmidt}), a group has modular subgroup lattice if and only if it is a direct product of groups of pairwise coprime orders each of which is either a modular $p$-group or a $P^*$-group. The modular $p$-groups are exactly the groups of types~(i) and~(ii) \cite[Theorem~2.3.1]{Schmidt} (see also \cite{Napolitani}). A $P^*$-group is a semidirect product of an elementary abelian $p$-group $V$ by a cyclic $q$-group $\gen{t}$, where $q\ne p$ and $t$ induces on $V$ a power automorphism of prime order. Every power automorphism of $V$ is the multiplication by some element of $\mathbb F_p^\times$: indeed, it leaves invariant every one-dimensional subspace of $V$, and considering the subspaces spanned by the sums of basis vectors we see that all the corresponding scalars coincide. Hence the $P^*$-groups are exactly the groups of type~(iii). The theorem is proved.
\end{proof}

\begin{corollary}\label{cor:fixed-commutator}
In Theorem~\ref{thm:main}, it is enough to consider only the subgroups $R$ such that $R=[P_R,D]$.
\end{corollary}

\begin{proof}
Let $R$ be as in Theorem~\ref{thm:main} and $W=[P_R,D]$. Then $W$ is normal in $D$, since $P_R$ is $D$-invariant, and $W\le R$. It is clear that $P_W=P_R$, and so $W=[P_W,D]$. Since $P_R/R$ is a quotient of $P_W/W$ and the subgroup lattice of a quotient of a group $X$ is isomorphic to an interval of $\Lall(X)$, the modularity of $\Lall(P_W/W)$ implies the modularity of $\Lall(P_R/R)$. The corollary is proved.
\end{proof}

\section{\texorpdfstring{Proof of Theorem~\ref{thm:nilpotent} and normal Hall subgroups}{Proof of Theorem B and normal Hall subgroups}}

\begin{lemma}\label{lem:nilpotent-invariant}
Let a group $D$ act on a nilpotent group $P$ by automorphisms. Then the lattice of all $D$-invariant subgroups of $P$ is modular if and only if every two of its members are permutable.
\end{lemma}

\begin{proof}
Let $\Lall$ be the lattice of all $D$-invariant subgroups of $P$. If every two members of $\Lall$ are permutable, then the join of any two members $A$ and $B$ of $\Lall$ is $AB$, and the modular law in $\Lall$ follows from the Dedekind identity.

Conversely, suppose that $\Lall$ is modular. First we show that if $U,V\in\Lall$ and $V$ covers $U$ in $\Lall$, then $U$ is normal in $V$. Indeed, $N_V(U)\in\Lall$ and, since $V$ is nilpotent, $U<N_V(U)$. Hence $N_V(U)=V$.

Now let $A,B\in\Lall$, and let $A\cap B=B_0<B_1<\cdots<B_n=B$ be a maximal chain in $\Lall$ from $A\cap B$ to $B$. Since $\Lall$ is modular, the map $X\mapsto\gen{A,X}$ is an isomorphism of the interval $[A\cap B,B]$ onto the interval $[A,\gen{A,B}]$ of $\Lall$, with the inverse map $Y\mapsto Y\cap B$. Hence the subgroups $A_k=\gen{A,B_k}$, $k=0,1,\ldots,n$, form a maximal chain in $\Lall$, and so $A_{k-1}$ is normal in $A_k$ for $k=1,\ldots,n$. It follows that $A_k=A_{k-1}B_k$, and by induction we have $A_k=AB_k$ for all $k$. In particular, $\gen{A,B}=A_n=AB$, and so $AB=BA$. The lemma is proved.
\end{proof}

\begin{proof}[Proof of Theorem~\ref{thm:nilpotent}]
If every two members of $\Lper(G)$ are permutable, then $\Lper(G)$ is modular by the Dedekind identity.

Now assume that $\Lper(G)$ is modular, and let $A,B\in\Lper(G)$, $N=A_GB_G$ and $Q=G/N$. Then $AN/N$ and $BN/N$ are $\sigma$-nilpotent, since they are quotients of $A/A_G$ and $B/B_G$, respectively, which are $\sigma$-nilpotent by Lemma~\ref{lem:skiba}(2). Moreover, $AN/N,BN/N\in\Lper(Q)$ by Lemma~\ref{lem:skiba}(4). Hence the Hall $\sigma_i$-subgroups of $AN/N$ and $BN/N$ belong to $\Lsi(Q)$ by Lemma~\ref{lem:skiba}(3). By Proposition~\ref{prop:quotient} and Lemma~\ref{lem:block}, $\Lsi(Q)$ is a modular lattice of $O^{\sigma_i}(Q)$-invariant subgroups of $O_{\sigma_i}(Q)$, and $O_{\sigma_i}(Q)$ is nilpotent by hypothesis. Therefore every two members of $\Lsi(Q)$ are permutable by Lemma~\ref{lem:nilpotent-invariant}. Since the subgroups $O_{\sigma_i}(Q)$ and $O_{\sigma_j}(Q)$ centralize each other for $i\ne j$, it follows that $AN/N$ and $BN/N$ are permutable, that is, $ABN=BAN$. But $A_G$ and $B_G$ are normal subgroups of $G$ contained in $A$ and $B$, respectively, so $ABN=AB$ and $BAN=BA$. Hence $AB=BA$.

Finally, if $H$ is a Hall $\sigma_i$-subgroup of $G$, then $HN/N$ is a Hall $\sigma_i$-subgroup of $G/N$, and $O_{\sigma_i}(G/N)\le HN/N$ (see the proof of Lemma~\ref{lem:block}). Hence if $H$ is nilpotent, then so is $O_{\sigma_i}(G/N)$. The theorem is proved.
\end{proof}

\begin{proposition}\label{prop:hall}
Let $G$ be a $\sigma$-full group and $V$ a normal Hall subgroup of $G$ such that $\sigma(V)\cap\sigma(G/V)=\varnothing$. If the lattices $\Lper(V)$ and $\Lper(G/V)$ are modular, then the lattice $\Lper(G)$ is modular. In particular, if $V$ is abelian, then $\Lper(G)$ is modular if and only if $\Lper(G/V)$ is modular.
\end{proposition}

\begin{proof}
If $\sigma_i\in\sigma(V)$, then every Hall $\sigma_i$-subgroup of $G$ is contained in $V$, since $\sigma_i\cap\pi(G/V)=\varnothing$. Hence $V$ is $\sigma$-full. Let $N$ be a normal subgroup of $G$, $Q=G/N$, $W=VN/N$ and $i\in I$. In view of Proposition~\ref{prop:quotient}, it is enough to show that the lattice $\Lsi(Q)$ is modular.

First assume that $\sigma_i\in\sigma(V)$. Since $Q/W$ is a $\sigma_i'$-group, we have $O_{\sigma_i}(Q)\le W$, and so $O_{\sigma_i}(Q)=O_{\sigma_i}(W)$. Moreover, $W/(W\cap O^{\sigma_i}(Q))$ is a $\sigma_i$-group, so $O^{\sigma_i}(W)\le O^{\sigma_i}(Q)$. Hence $\Lsi(Q)$ is a sublattice of $\Lsi(W)$ by Lemma~\ref{lem:block}. Since $W\simeq V/(V\cap N)$ and the lattice $\Lper(V)$ is modular, the lattice $\Lsi(W)$ is modular by Proposition~\ref{prop:quotient}. Therefore $\Lsi(Q)$ is modular.

Now assume that $\sigma_i\notin\sigma(V)$. Then $W$ is a $\sigma_i'$-group, so $O_{\sigma_i}(Q)\cap W=1$ and $W\le O^{\sigma_i}(Q)$. Hence the natural homomorphism $Q\to Q/W$ maps $O_{\sigma_i}(Q)$ isomorphically onto a normal $\sigma_i$-subgroup $\overline P$ of $Q/W$, and it maps $O^{\sigma_i}(Q)$ onto $O^{\sigma_i}(Q/W)$. Moreover, a subgroup $U$ of $O_{\sigma_i}(Q)$ is $O^{\sigma_i}(Q)$-invariant if and only if $UW/W$ is $O^{\sigma_i}(Q/W)$-invariant, since $UW\cap O_{\sigma_i}(Q)=U$. Therefore $\Lsi(Q)$ is isomorphic to the interval $[1,\overline P]$ of $\Lsi(Q/W)$ by Lemma~\ref{lem:block}. Since $Q/W\simeq G/VN$ is a quotient of $G/V$ and the lattice $\Lper(G/V)$ is modular, the lattice $\Lsi(Q/W)$ is modular by Proposition~\ref{prop:quotient}. Therefore $\Lsi(Q)$ is modular.

Hence $\Lper(G)$ is modular. Finally, if $V$ is abelian, then the lattice $\Lper(V)$ is modular, being a sublattice of the modular lattice $\Lall(V)$; and if $\Lper(G)$ is modular, then so is $\Lper(G/V)$ by Lemma~\ref{lem:skiba}(4). The proposition is proved.
\end{proof}

\Needspace{9\baselineskip}
\section{Examples}

\begin{example}\label{ex:dihedral}
Let $\sigma$ be a partition of $\mathbb P$ such that $2$ and $3$ belong to the same member $\sigma_1$ of $\sigma$ and $5\notin\sigma_1$. Let
\[
G=D_{30}=\gen{r,s\mid r^{15}=s^2=1,\ srs=r^{-1}}.
\]
Then the lattice $\Lper(G)$ is modular, but not every two members of $\Lper(G)$ are permutable.

Indeed, $V=\gen{r^3}\simeq C_5$ is a normal Hall subgroup of $G$ and $G/V\simeq S_3$ is a $\sigma_1$-group, so $\sigma(V)\cap\sigma(G/V)=\varnothing$ and $\Lper(G/V)=\Lall(G/V)$. The lattice $\Lall(S_3)$ has length $2$, so it is modular. Hence $\Lper(G)$ is modular by Proposition~\ref{prop:hall}. The subgroups $V\gen{s}$ and $V\gen{r^5s}$ of order $10$ belong to $\Lper(G)$ by Lemma~\ref{lem:skiba}(4), since their images in $G/V$ belong to $\Lper(G/V)=\Lall(G/V)$. But $V\gen{s}\cap V\gen{r^5s}=V$, so the product of these subgroups has $20$ elements, and hence these subgroups are not permutable. Note that the Hall $\sigma_1$-subgroups of $G$ are isomorphic to $S_3$, so they are not nilpotent.
\end{example}

\begin{example}\label{ex:quotients}
The lattices $\Lsi(G)$, $i\in I$, may be modular, while the lattice $\Lper(G)$ is not modular. Thus in Proposition~\ref{prop:quotient} it is not enough to consider only the lattices $\Lsi(G)$.

Indeed, let $V=\mathbb F_3^2$ and $Q=\gen{r,s}\le\mathrm{GL}_2(3)$, where
\[
r=\begin{pmatrix}0&-1\\1&0\end{pmatrix}\qquad\text{and}\qquad
s=\begin{pmatrix}1&0\\0&-1\end{pmatrix}.
\]
Then $Q\simeq D_8$, and $Q$ acts faithfully and irreducibly on $V$, since the characteristic polynomial $x^2+1$ of $r$ is irreducible over $\mathbb F_3$. Let $G=V\rtimes Q$ and $\sigma=\{\{2\},\{3\},\{5\},\ldots\}$. If $X$ is a normal $2$-subgroup of $G$, then $[X,V]\le X\cap V=1$, so $X\le C_G(V)=V$ and hence $X=1$. Thus $O_2(G)=1$. Moreover, $O_3(G)=V$ and $Q\le O^3(G)$. Hence $\Lall_{\{2\}}(G)=\{1\}$ and $\Lall_{\{3\}}(G)=\{1,V\}$, and $\Lall_{\{p\}}(G)=\{1\}$ for every prime $p>3$, so all these lattices are modular. On the other hand, by Lemma~\ref{lem:skiba}(4), the interval of $\Lper(G)$ consisting of all members containing $V$ is isomorphic to $\Lper(G/V)=\Lall(D_8)$, which is not modular. Indeed, in $D_8=\gen{r,s}$ let $a=\gen{s}$, $b=\gen{rs}$ and $c=\gen{s,r^2}$. Then $a\le c$ and $b\cap c=1$, but $\gen{a,b}\cap c=c\ne a=\gen{a,b\cap c}$.
\end{example}

\begin{example}\label{ex:intermediate}
In Theorem~\ref{thm:main}, it is not enough to consider only the subgroups $R=1$ and $R=P\cap D$.

Indeed, let $H=D_8$, $P=H\times H\times H$ and $G=P\rtimes\gen{c}$, where $c$ has order $3$ and permutes the factors of $P$ cyclically. Let $\sigma=\{\{2\},\{3\},\{5\},\ldots\}$ and $\sigma_i=\{2\}$. Then $O_2(G)=P$ and $O^2(G)=D=K\rtimes\gen{c}$, where
\[
K=\{(h_1,h_2,h_3)\in P\mid h_1h_2h_3\in H'\}.
\]
Indeed, the map $(h_1,h_2,h_3)c^j\mapsto h_1h_2h_3H'$ is a homomorphism of $G$ onto the $2$-group $H/H'$ with kernel $D$, so $O^2(G)\le D$. On the other hand, in every $2$-group quotient of $G$ the image of $c$ is trivial, so the three factors of $P$ have the same image, which is abelian since these factors centralize each other; hence $D\le O^2(G)$. Note that $P\cap D=K$.

For $R=K$, we have $P_R=P$, since $[P,D]\le P\cap D=K$, and so $P_R/R=P/K\simeq C_2\times C_2$. For $R=1$, we have $P_R=C_P(D)$. Every element of $C_P(D)$ centralizes $c$, so it has the form $(h,h,h)$; moreover, $h\in Z(H)$, since $(u,u^{-1},1)\in K$ for all $u\in H$. Hence $C_P(D)\simeq C_2$. Thus in both cases the group $P_R/R$ has modular subgroup lattice.

Now let $R=\{(z_1,z_2,z_3)\in Z(H)^3\mid z_1z_2z_3=1\}$. Then $R$ is a normal subgroup of $G$ and $R\le K$. Let $h\in H$ and $k=(u_1,u_2,u_3)\in K$. The commutator of $(h,h,h)$ and $k$ is $([h,u_1],[h,u_2],[h,u_3])$, its coordinates lie in $H'=Z(H)$, and their product is $[h,u_1u_2u_3]=1$, since $u_1u_2u_3\in H'=Z(H)$. Moreover, $(h,h,h)$ centralizes $c$. Hence the diagonal subgroup $\Delta=\{(h,h,h)\mid h\in H\}\simeq D_8$ is contained in $P_R$. Since $|Z(H)|=2$ and $z^3=z$ for all $z\in Z(H)$, we have $\Delta\cap R=1$. Therefore $P_R/R$ has a subgroup isomorphic to $D_8$, and so its subgroup lattice is not modular.
\end{example}

\end{document}